\documentclass{ifacconf}

\usepackage{graphicx}      
\usepackage{amsmath}
\usepackage{bbm}
\usepackage{epstopdf}
\usepackage{tikz}
\usepackage{natbib}
\usetikzlibrary{shapes,positioning,backgrounds}
\RequirePackage{amsmath}
\RequirePackage{amssymb}
\begin{document}
\begin{frontmatter}

\title{Stochastic optimal control for nonlinear damped network dynamics} 


\author[First]{Simone G\"ottlich} 
\author[Second]{Thomas Schillinger} 

\address[First]{University of Mannheim, Department of Mathematics, 68131 Mannheim, Germany (e-mail: goettlich@uni-mannheim.de)}
\address[Second]{University of Mannheim, Department of Mathematics, 68131 Mannheim, Germany (e-mail: schillinger@uni-mannheim.de)}

\begin{abstract}                
We present a stochastic optimal control problem for a tree network. The dynamics of the network are governed by transport equations with a special emphasis on the non-linear damping function. Demand profiles at the network sinks are modelled by a stochastic differential equations. An explicit optimal inflow into the network is determined and numerical simulations are presented to show the effects for different choices of the non-linear damping.
\end{abstract}


\end{frontmatter}

\tikzset{
	es/.style={
		circle,
		fill=white,
		draw=black
	}
}

\section{Introduction}
Energy networks or supply networks in a more general setting can be analyzed from different perspectives. Mathematically, there is a strong interest to describe the dynamics inside the network to better understand the underlying physical or economic processes. Typical applications range from electric transmission lines \cite{Goettlich2015} and gas networks \cite{Banda2006} to production systems \cite{Apice}, wherein the dynamics is governed by nonlinear hyperbolic transport equations \cite{Bressan2014}. Generally, nonlinearities might appear due to complex flow patterns or additional interaction terms such as resistance or friction.     
For many applications the question arises how such a network can be controlled to satisfy consumer demands. Since demands usually include a kind of uncertainty, this leads to a challenging stochastic optimal control problem. To reduce the complexity of the problem, we therefore start with the consideration of linearized dynamics and nonlinear damping combined with stochasticity of demands in an abstract setting. This enables us to investigate the qualitative behaviour of
solutions to the optimal control problem.    

This work is based on ideas originally presented in \cite{GoettlichSchillinger,Lux1}. Basically, we consider a network framework consisting of three ingredients. At the network sources an optimal inflow shall be injected into the system such that a given demand is met. On every arc of the network transport equations of hyperbolic type describe the dynamics. In our case we focus on linear flux functions and add unlike \cite{GoettlichSchillinger} a non-linear damping term. However, we are still able to compute an explicit representation of the optimal input. 
As a last component, we consider uncertain demands at the sinks of the network. These are described by a stochastic process, given by a stochastic differential equation. In this work, we assume that demands are described by Jacobi processes, as recently proposed by \cite{Korn2021}. These are mean-reverting processes which stay in a bounded interval and are therefore very suitable for various applications we have in mind. The goal of this contribution is to show under which  assumptions we will be able to explicitly derive the optimal inflow subject to the network dynamics and the stochastic demand.

The organization of the article is the following. In Section \ref{secOptCont} we present the full optimal control framework. The discussion of the objective function and the availability of information is executed in Section \ref{solObj}. The key result is presented in Section \ref{nonlinearDamping}, where we extensively describe how an explicit formula for the optimal inflow can be calculated in the setting of non-linear damping terms. Section \ref{DPa} deals with the choices of  distribution parameters at intersections and is followed by a numerical study of different non-linear damping functions in Section \ref{numSim}.

\section{The Optimal Control Problem}
\label{secOptCont}
We consider the control problem in (\ref{ControlProblem11}), where we aim to determine the optimal input $u(t)$ given some stochastic demand $D_s^{(v_i)}$ and the network dynamics in terms of scalar transport equations.
More precisely, we make use of the following notation: For a directed tree network, the set $C$ denotes the demand nodes and $J$ the set of the (inner) intersections. We control the network inflow $u(t)$ defined in (\ref{boundaryData}) and solve an optimization problem for the demand nodes. At any demand node $v_i \in C$ we consider a demand process $(D_t^{(v_i)})_{t \in [t_0,T]}$ for which we assume a Jacobi process. The Jacobi process can be described by a stochastic differential equation given in (\ref{Demand})-(\ref{initialDemand}), where $\theta^{(v_i)}$ is a time-dependent mean reversion level, $\kappa^{(v_i)}$ the mean reversion speed and $\sigma^{(v_i)}$ the scaling of the influence of stochastic disturbances coming from a Brownian motion $(W_t^{(v_i)})_{t \in [t_0,T]}$. In this work, we only consider Jacobi processes on $[0,1]$, which can be easily translated on any bounded interval. For further details of the Jacobi process see e.g. \cite{Sorensen} or \cite{JacobiG}. 
\begin{subequations}\label{ControlProblem11}
	\begin{align}
		\underset{u \in L^2}{\min} & \sum_{\lbrace i: v_i \in C \rbrace} \int_{t_0}^T \mathbb{E}\left[\left(D_s^{(v_i)} - f^{(i)}(z^{(i)}(1,s),s)\right)^2 ~\Big|~ \mathcal{F}_{\hat{t}}\right] ds  \label{ObjectiveFunc}\\
		\text{s.t.}~~ & z_t^{(i)}(x,t) + f^{(i)}(z^{(i)}(x,t),t)_x + g^{(i)}(z^{(i)}(x,t),t) = 0 \label{NetworkDynamics}\\ 
		& z^{(i)}(x,t_0) = z_0^{(i)}(x), ~~~~ \forall i \text{ s.t. } v_i \in J \cup C \label{initialData}\\
		& f^{(1)}(z^{(1)}(0,t),t) = u(t) \label{boundaryData}\\
		&f^{(k)}(z^{(k)}(0,t),t) =  \alpha_{i,k}(t)f^{(i)}(z^{(i)}(1,t),t), \nonumber\\
		&~~ \forall i \text{ s.t. } v_i \in J,~ k\text{ outgoing arc of }v_i \label{fluxCons1}\\
		& \sum_{k \text{ outgoing arc of } v_i} \alpha_{i,k}(t) = 1,~~~ \forall i \text{ s.t. } v_i \in J\label{fluxCons2}\\
		& dD_t^{(v_i)} = \kappa^{(v_i)} \left(\theta^{(v_i)}(t) - D_t^{(v_i)}\right) \nonumber\\
		&~~~+ \sigma^{(v_i)} \sqrt{D_t^{(v_i)}\left(1-D_t^{(v_i)}\right)} dW_t^{(v_i)} \label{Demand}\\
		& ~D_0^{(v_i)} = d_0^{(v_i)},~ v_i \in C . \label{initialDemand}
	\end{align}
\end{subequations}

We aim to minimize the expected quadratic deviation between the demand process and the outflow of the network (\ref{ObjectiveFunc}). In the objective function, we condition on a time $\hat{t} \leq t_0$ up to which demand information is available. Section \ref{solObj} presents further strategies about demand updates.  On the network arcs $i$ we consider hyperbolic partial differential equations for the quantities $z^{(i)}(x,t)$ in (\ref{NetworkDynamics}), where the flux functions $f^{(i)}$ are linear with respect to $z^{(i)}(x,t)$ and possibly non-linear in time, accompanied with some initial data (\ref{initialData}). In comparison to the work in \cite{GoettlichSchillinger}, the damping functions $g^{(i)}$ are also allowed the be non-linear in the quantities $z^{(i)}$. The choices of the damping functions will be discussed in Section \ref{nonlinearDamping}. At each inner node due to the tree network structure, there is exactly one ingoing arc but possibly several outgoing arcs. To ensure flux conservation, we have to impose distributions parameters $\alpha_{i,k}(t)$ for the share of flux moving from the ingoing arc $i$ to an outgoing arc $k$ in (\ref{fluxCons1})-(\ref{fluxCons2}). The distribution parameters are not the main focus of this work but will shortly be discussed in Section \ref{DPa}.

\section{Solution to the Objective Function and Demand Updates}
\label{solObj}
The optimal control problem in (\ref{ControlProblem11}) has a stochastic component given by the stochastic differential equation (SDE) for the demand process (\ref{Demand}). But when minimizing the expected quadratic deviation of demand and network supply the optimal solution for a square integrable demand process is given by $\mathbb{E}[D_t^{(v_i)} | \mathcal{F}_{\hat{t}} ]$ (see Corollary 8.17 in \cite{Klenke}). Therefore, the particular structure of the demand process, apart from the conditional expectation, does not matter to solve the optimal control problem. Considering the Jacobi process the conditional expectation with respect to $t_0$ is explicitly given by (see e.g. \cite{Delbaen2002})
\begin{align}
\begin{split}
	\label{condfirstMomentTimevar}
	\mathbb{E}[D_t^{(v_i)} | D_{t_0}^{(v_i)} = d_0^{(v_i)}] &= d_0^{(v_i)} e^{-\kappa^{(v_i)}(t-t_0)} \\
	&~~~+ \kappa^{(v_i)} \int_{t_0}^t e^{-\kappa^{(v_i)}(s-t_0)}\theta^{(v_i)}(s)ds.
\end{split}
\end{align}
As a next step, we want to generalize the control problem a little further. So far, we have assumed that there is one time $\hat{t}\leq t_0$ at which the demand levels are updated. However, we now consider a sequence of update times $(\hat{t}_j)_{j \in \mathbb{N}}$ which allow for additional information about the demands. Then, we end up with a sequence of optimal control problems on time intervals $[\hat{t}_j,\hat{t}_{j+1})$. To account for the correct time intervals, the objective function from (\ref{ObjectiveFunc}) then reads
\begin{align*}
	\underset{u \in L^2}{\min}  \sum_{\lbrace i: v_i \in C \rbrace} &\int_{\tilde{t}(v_i,\hat{t}_j)}^{\tilde{t}(v_i,\hat{t}_{j+1})} \mathbb{E}\Bigg[\bigg(D_s^{(v_i)} \\
	&- f^{(i)}(z^{(i)}(1,s),s)\bigg)^2 ~\Big|~ \mathcal{F}_{\hat{t}_j}\Bigg] ds ,
\end{align*}
where $\tilde{t}(v_i,t)$ denotes the time at which the inflow inserted at the source node at time $t$ reaches the demand node $v_i$. Additionally, for all updates except the first one, the initial data has to correspond to the state of the system of the time interval before. Therefore, the equations (\ref{initialData}) and (\ref{initialDemand}) read
\begin{align*}
	z^{(i)}(x,\hat{t}_j) &= z_\text{old}^{(i)}(x,\hat{t}_j), ~~~~ \forall i \text{ s.t. } v_i \in J \cup C \\
	D_{\hat{t}_j}^{(v_i)} &= D_{\hat{t}_j}^{(v_i),\text{old}},~ v_i \in C ,
\end{align*}
where $z_\text{old}^{(i)}(x,\hat{t}_j)$ and $D_{\hat{t}_j}^{(v_i),\text{old}}$ denote the final quantity of arc $i$ at position $x$ and the final demand at node $v_i$ in the previous time interval, respectively.

\section{Network dynamics with non-linear damping}
\label{nonlinearDamping}
In this section, we focus on the dynamics in the network on all arcs $i$ given by constraint (\ref{NetworkDynamics}) of the optimal control problem, i.e. the shape of the functions $f^{(i)}$ and $g^{(i)}$ governed by
\begin{align}
	\label{PDEwD}
	z_t^{(i)}(x,t) + f^{(i)}(z^{(i)}(x,t),t)_x + g^{(i)}(z^{(i)}(x,t),t) = 0,
\end{align}
where $f^{(i)}$ denotes the flux function and $g^{(i)}$ is the damping function. 
Here, we restrict to functions $f^{(i)}$ that are linear in $z^{(i)}$ but non-linear in  the time and potentially non-linear functions $g^{(i)}$ in $z^{(i)}$ and $t$. 
The damping reflects a loss in the transported quantity over time, which may be due to some physical property as for instance friction or electrical resistance. 
The linearity of the flux-functions will play a very important role because characteristic curves do not cross and hence no discontinuities may appear, see \cite{LeVeque} for an overview. For simplicity, all investigations are executed in the 1-1 and 1-2-network case and can be generalized to arbitrary tree networks in a straightforward way. 

In the following, we discuss the influence of the non-linear damping term for the two relevant cases. We restrict to non-linear damping functions in which the non-linearities of time and quantity are separated, i.e. $$g^{(i)}(z,t) = \mu_i(t) \hat{g}_i(z),$$ 
where $\mu_i \in L^1([t_0,T])$ is chosen to be a non-negative function and $\hat{g}_i$ is a Lipschitz-continuous function whose antiderivative is explicitly known. Additionally, we assume $\hat{g}_i\geq 0$ and $\hat{g}_i=0$ only on a Lebesgue null-set. This allows for the natural choice of $\hat{g}_i(0)=0$, but still preserves some important properties we will exploit. The flux function on arc $i$ is chosen as  $f^{(i)}(z,t) = \lambda_i(t) z$ with strictly positive velocity $\lambda_i(\cdot)$ . 

\subsubsection{1-1 Network}
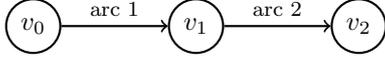
\begin{figure}[ht!]
	\label{12network}
	\begin{center}
		\begin{tikzpicture}[thick]
			\node[es] (K1) {$v_0$};
			\node[es,right = 1.4cm of K1] (K2) {$v_1$};
			\node[es,right = 1.4cm of K2] (K3) {$v_2$};

			\path [->]
			(K1) edge node[above]{\small{arc 1}} (K2)
			(K2) edge node[above]{\small{arc 2}}(K3);
			
		\end{tikzpicture}
	\end{center}
	\caption{The supply system as a 1-1-network with one source $v_0$ and one demand node $v_2$.}
	\label{network1-1}
\end{figure}

First, we show how the optimal inflow for a 1-1 network has to be chosen if we consider a possible non-linear damping term on both arcs. Denote by $\mu_i(t)$ a time-dependent damping factor for arc $i$ and $\hat{g}_i(z) = \frac{1}{\tilde{g}_i(z)}$ the non-linear damping function for arc $i$ such that the overall damping term is given by $g_i(z,t) = \frac{\mu_i(t)}{\tilde{g}_i(z)}$. We make use of the network from Figure $\ref{network1-1}$, where we have one source node $v_0$, an inner intersection $v_1$ and a demand node $v_2$. 

\label{1-1damp}
To relate the injection time at the source node $v_0$ and the output time at the demand node $v_2$ we consider a characteristic curve $(x(t),t)$ of a unit for which $t_2$ is the time in which the unit injected at $t_{\text{in}}$ reaches the demand node $v_2$ and $t_1$ the time when it reaches node $v_1$. If the velocity functions $\lambda_i$ have an antiderivative $\Lambda_i, i =1,2$ which is invertible, for some injection time $t_\text{in}$ the values of $t_1$ and $t_2$ are explicitly given by
\begin{align}
	\begin{split}
		\label{intermediateTimes2}
		t_1 &= \Lambda_1^{-1}(1+\Lambda_1(t_{\text{in}}))\\
		t_2 &= \Lambda_2^{-1}(1+\Lambda_2(\Lambda_1^{-1}(1+\Lambda_1(t_{\text{in}})))),
	\end{split}
\end{align}
where we for simplicity assume that every arc has length 1 (see \cite{GoettlichSchillinger}).


Staying on such a characteristic curve, in a setting without damping, the transported quantity $z(x(t),t)$ remains constant. If we introduce the damping term, the quantity $z$ is reduced by the damping. This can be formulated by the following ordinary differential equation (ODE)
\begin{align}
	\label{ODEdamp}
	\frac{d}{dt}z^{(i)}(x(t),t) = -g_i(z^{(i)}(x(t),t),t) = -\frac{\mu_i(t)}{\tilde{g}_i(z^{(i)}(x(t),t))}.
\end{align}
This ODE can be uniquely solved due the Lipschitz-continuity of $\tilde{g}$. The technique of separation of variables even allows for the calculation of an explicit solution. To do so, we assume that an antiderivative of $\tilde{g}_i$ is explicitly given by $\tilde{G}_i$ and we consider time points $t_{k-1} < t_k$:
\begin{align*}
	&\frac{d}{dt}z^{(i)}(x(t),t) = -g_i(z^{(i)}(x(t),t),t) = -\frac{\mu_i(t)}{\tilde{g}_i(z^{(i)}(x(t),t))}\\
	\Leftrightarrow &~\frac{d}{dt}\left(z^{(i)}(x(t),t)\right) \cdot \tilde{g}_i(z^{(i)}(x(t),t)) = -\mu_i(t).
\end{align*}
Integrating both sides of the equation over $t_{k-1}$ to $t_k$ and applying the fundamental theorem of calculus leads to
\begin{align*}
	&\int_{t_{k-1}}^{t_k} \frac{d}{ds}\left(z^{(i)}(x(s),s)\right) \cdot \tilde{g}_i(z^{(i)}(x(s),s))ds \\
	&= -\int_{t_{k-1}}^{t_k} \mu_i(s) ds\\
	\Leftrightarrow &~ \left[ \tilde{G}_i\left(z^{(i)}(x(s),s)\right) \right]_{t_{k-1}}^{t_k} = -\int_{t_{k-1}}^{t_k} \mu_i(s) ds\\
	\Leftrightarrow & ~\tilde{G}_i\left(z^{(i)}(x(t_{k-1}),t_{k-1})\right) \\
	&= \tilde{G}_i\left(z^{(i)}(x(t_{k}),t_{k})\right) + \int_{t_{k-1}}^{t_k} \mu_i(s) ds\\
	\Leftrightarrow &~ z^{(i)}(x(t_{k-1}),t_{k-1}) \\
	&= \tilde{G}_i^{-1} \left(\tilde{G}_i\left(z^{(i)}(x(t_{k}),t_{k})\right) + \int_{t_{k-1}}^{t_k} \mu_i(s) ds \right).
\end{align*}
The antiderivative of $\tilde{g}_i$ is continuous and strictly increasing because $\tilde{g}_i$ is strictly positive, apart from a Lebesgue null-set. Then, there exists the inverse function $(\tilde{G}_i)^{-1}$ of $\tilde{G}_i$ and the solution is well-defined. The integral of the damping function stays bounded since $\mu_i \in L^1([t_0,T])$. Note that we are going to perform a backward calculation for the optimal inflow under the assumption that we know the corresponding optimal outflow of the network. Therefore, we consider an end value problem here, this means we have upstream information at time $t_k$ and calculate the corresponding initial values at time $t_{k-1}$. 
We now present this backward calculation for the optimal inflow and start at the demand node $v_2$.
\begin{align*}
	&z^{(2)}(1,t_2) \lambda_2(t_2) \Leftrightarrow z^{(2)}(1,t_2) = \frac{f^{(2)}(z^{(2)}(1,t_2),t_2)}{\lambda_2(t_2)}
\end{align*}
Following the characteristic curve and using the result of the ODE solution we can calculate the resulting quantity at the beginning of arc 2 at time $t_1$ by
\begin{align*}
	&z^{(2)}(0,t_1) = \tilde{G}_2^{-1} \left(\tilde{G}_2\left(z^{(2)}(1,t_{2})\right) + \int_{t_{1}}^{t_2} \mu_2(s) ds \right) \\
	&= \tilde{G}_2^{-1} \left(\tilde{G}_2\left(\frac{f^{(2)}(z^{(2)}(1,t_2),t_2)}{\lambda_2(t_2)}\right) + \int_{t_{1}}^{t_2} \mu_2(s) ds \right).
\end{align*}
Using this, we can calculate the ingoing flux into arc 2 at time $t_1$ by
\begin{align*}
	&f^{(2)}(z^{(2)}(0,t_1),t_1) = \lambda_2(t_1) z^{(2)}(0,t_1)  \\
	&= \lambda_2(t_1)\tilde{G}_2^{-1} \left(\tilde{G}_2\left(\frac{f^{(2)}(z^{(2)}(1,t_2),t_2)}{\lambda_2(t_2)}\right) + \int_{t_{1}}^{t_2} \mu_2(s) ds \right).
\end{align*}
Accounting for flux conservation at the intersection it must hold that 
\begin{align*}
	f^{(2)}(z^{(2)}(0,t_1),t_1) = f^{(1)}(z^{(1)}(1,t_1),t_1).
\end{align*}
Then, we can deduce the corresponding quantity at the end of arc 1 by
\begin{align*}
	z^{(1)}(1,t_1) &= \frac{f^{(1)}(z^{(1)}(1,t_1),t_1)}{\lambda_1(t_1)} \\
	&= \frac{\lambda_2(t_1)}{\lambda_1(t_1)}\cdot \tilde{G}_2^{-1} \bigg(\tilde{G}_2\left(\frac{f^{(2)}(z^{(2)}(1,t_2),t_2)}{\lambda_2(t_2)} \right) \\
	&~~+ \int_{t_{1}}^{t_2} \mu_2(s) ds \bigg).
\end{align*}
Now again applying the ODE solution on arc 1 we get for the initial quantity
\begin{align*}
	&z^{(1)}(0,t_{\text{in}}) \\
	&= \tilde{G}_1^{-1} \bigg(\tilde{G}_1 \bigg[ \frac{\lambda_2(t_1)}{\lambda_1(t_1)}\cdot\tilde{G}_2^{-1}\bigg(\tilde{G}_2\left(\frac{f^{(2)}(z^{(2)}(1,t_2),t_2)}{\lambda_2(t_2)} \right) \\
	&~~+ \int_{t_{1}}^{t_2} \mu_2(s) ds \bigg)\bigg]  + \int_{t_{\text{in}}}^{t_1} \mu_1(s)ds \bigg).
\end{align*}
The ingoing flux into the 1-1 network is now a direct consequence of the quantity at the beginning of arc 1 and given by
\begin{align*}
	&f^{(1)}(z^{(1)}(0,t_{\text{in}}),t_{\text{in}}) \\
	&= \lambda_1(t_{\text{in}}) z^{(1)}(0,t_{\text{in}}) \\
	&= \lambda_1(t_{\text{in}})\tilde{G}_1^{-1} \bigg(\tilde{G}_1 \bigg[ \frac{\lambda_2(t_1)}{\lambda_1(t_1)} \tilde{G}_2^{-1} \Big(\tilde{G}_2\left(\frac{f^{(2)}(z^{(2)}(1,t_2),t_2)}{\lambda_2(t_2)} \right) \\
	&~~~~+ \int_{t_{1}}^{t_2} \mu_2(s) ds \Big)\bigg]  + \int_{t_{\text{in}}}^{t_1} \mu_1(s)ds \bigg).
\end{align*}

Since the optimal outflow of the network should match the corresponding conditional expected demand we obtain for the optimal inflow 
\begin{align*}
	&u(t_{\text{in}}) \\
	&= \lambda_1(t_{\text{in}})\cdot \tilde{G}_1^{-1} \Bigg(\tilde{G}_1 \bigg[ \frac{\lambda_2(t_1)}{\lambda_1(t_1)}\cdot \tilde{G}_2^{-1} \bigg(\tilde{G}_2\left(\frac{\mathbb{E}\left[D_{t_2}^{(v_2)} | \mathcal{F}_{\hat{t}}\right]}{\lambda_2(t_2)} \right) \\
	&~~~~+ \int_{t_{1}}^{t_2} \mu_2(s) ds \bigg)\bigg]  + \int_{t_{\text{in}}}^{t_1} \mu_1(s)ds \Bigg),
\end{align*}
where $\hat{t}<t_{\text{in}}$ is the time of the latest demand update. Iteratively, this procedure can be extended to larger 1-1 networks.

\subsubsection{1-2 Network}
\begin{figure}[h]
	\label{12network1}
	\begin{center}
		\begin{tikzpicture}[thick]
			\node[es] (K1) {$v_0$};
			\node[es,right = 0.8cm of K1] (K2) {$v_1$};
			\node[es,above right = 0.8cm of K2] (K3) {$v_2$};
			\node[es, below right = 0.8cm of K2] (K4) {$v_3$};

			\path [->]
			(K1) edge node[above]{\small{arc 1}} (K2)
			(K2) edge node[above left]{\small{arc 2}}(K3)
			(K2) edge node[above right]{\small{arc 3}}(K4);
		\end{tikzpicture}
	\end{center}
	\caption{The supply system as a 1-2-network with one source and two demand nodes.}
	\label{network1-2}
\end{figure}
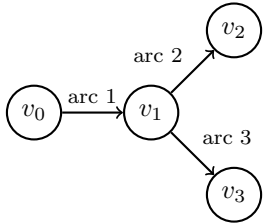

To expand the computations on arbitrary tree networks we exemplary present how to proceed in a 1-2 network as given in Figure \ref{network1-2} with one source node $v_0$, one intersection $v_1$ and two demand nodes $v_2$, $v_3$. 

As a direct consequence of the calculations for the 1-1-network, the ingoing fluxes into arc 2 and 3 at time $t_1$ have to amount to
\begin{align}
	\begin{split}
		\label{f2f3}
		f^{(2)}(0,t_1) &= \lambda_2(t_1) \cdot\tilde{G}_2^{-1} \bigg(\tilde{G}_2\left(\frac{f^{(2)}(z^{(2)}(1,t_2),t_2)}{\lambda_2(t_2)}\right) \\
		&~~+ \int_{t_1}^{t_2} \mu_2(s) ds  \bigg) \\
		f^{(3)}(0,t_1) &= \lambda_3(t_1) \cdot\tilde{G}_3^{-1} \bigg(\tilde{G}_3\left(\frac{f^{(3)}(z^{(3)}(1,t_3),t_3)}{\lambda_3(t_3)}\right) \\
		&~~+ \int_{t_1}^{t_3} \mu_3(s) ds  \bigg).
	\end{split}
\end{align}
As before $t_2$ is the time at which a unit injected at $t_\text{in}$ reaches node $v_2$ and similarly $t_3$ represents the time when it reaches node $v_3$.
To ensure flux conservation at the inner node we require 
\begin{align*}
	f^{(1)}(z^{(1)}(1,t_1)) = f^{(2)}(z^{(2)}(1,t_2)) + f^{(3)}(z^{(3)}(1,t_3)),
\end{align*}
which leads to an end quantity on arc 1 at time $t_1$ of
\begin{align*}
	z^{(1)}(1,t_1) &= \frac{f^{(1)}(z^{(1)}(1,t_1))}{\lambda_1(t_1)} \\
	&= \frac{f^{(2)}(z^{(2)}(1,t_2)) + f^{(3)}(z^{(3)}(1,t_3))}{\lambda_1(t_1)}.
\end{align*}
By the solution to the damped ODE (\ref{ODEdamp}) we get that the value of the quantity on the corresponding characteristic curve at time $t_{\text{in}}$ is given by
\begin{align*}
	&z^{(1)}(0,t_{\text{in}}) \\
	&= \tilde{G}_1^{-1} \bigg( \tilde{G}_1\left(\frac{f^{(2)}(z^{(2)}(1,t_2)) + f^{(3)}(z^{(3)}(1,t_3))}{\lambda_1(t_1)} \right) \\
	&~~~+ \int_{t_{\text{in}}}^{t_1} \mu_1(s)ds \bigg).
\end{align*}
As a last step we determine the corresponding ingoing flux and insert the equations for the ingoing fluxes into arc 2 and 3 from (\ref{f2f3}). Additionally using the information that the optimal outflow should match the conditional demand expectation we get an explicit form for the optimal inflow $u(t_{\text{in}})$ at time $t_{\text{in}}$ by
\begin{align*}
	&u(t_{\text{in}}) = \lambda_1(t_{\text{in}})z^{(1)}(0,t_{\text{in}}) \\
	& = \lambda_1(t_{\text{in}}) \cdot \tilde{G}_1^{-1} \bigg( \tilde{G}_1\left(\frac{f^{(2)}(z^{(2)}(1,t_2)) + f^{(3)}(z^{(3)}(1,t_3))}{\lambda_1(t_1)} \right) \\
	&~~+ \int_{t_{\text{in}}}^{t_1} \mu_1(s)ds \bigg)\\
	&= \lambda_1(t_{\text{in}}) \cdot \tilde{G}_1^{-1} \Bigg[ \tilde{G}_1\Bigg[\frac{\lambda_2(t_1)}{\lambda_1(t_1)} \cdot\tilde{G}_2^{-1} \bigg(\tilde{G}_2\left(\frac{\mathbb{E}\left[D_{t_2}^{(v_2)} |\mathcal{F}_{\hat{t}}\right]}{\lambda_2(t_2)}\right) \\
	&~~+ \int_{t_1}^{t_2} \mu_2(s) ds  \bigg) \\
	&~~+ \frac{\lambda_3(t_1)}{\lambda_1(t_1)} \cdot\tilde{G}_3^{-1} \left(\tilde{G}_3\left(\frac{\mathbb{E}\left[D_{t_3}^{(v_3)} |\mathcal{F}_{\hat{t}}\right]}{\lambda_3(t_3)}\right) + \int_{t_1}^{t_3} \mu_3(s) ds  \right) \Bigg] \\
	&~~+ \int_{t_{\text{in}}}^{t_1} \mu_1(s)ds \Bigg].
\end{align*}
This procedure can now be adapted consecutively also to larger tree networks. 

\section{Distribution Parameters}
\label{DPa}
At every inner node of the network we request flux conservation, meaning that all ingoing flux has to be distributed to the outgoing arcs, cf. Figure \ref{network1-2}. In the optimal control problem this is represented by the equations (\ref{fluxCons1}) and (\ref{fluxCons2}). The choice of the distribution parameters $\alpha_{i,k}(s)$ will be explained in this section. Assuming that we are in a setting with updated inflows, we set the distribution parameters for the flux proportional to the inflow shares for the corresponding demand nodes at the time $t$ at which the units where injected into the network. This can be understood as tracking the characteristic curve at node $v_i$ and time $s$ back to the source node. In the exemplary case of the 1-2 network and the notation from Section \ref{nonlinearDamping} the distribution parameters $\alpha_{1,k}(t)$ for $k=2,3$ are given by
\begin{align}
\label{distrPa}
	\alpha_{1,k} (t_1) = \frac{\lambda_k(t_1)\tilde{G}_k^{-1} \bigg(\tilde{G}_k\left(\frac{\mathbb{E}\left[D_{t_k}^{(v_k)} | \mathcal{F}_{\hat{t}_j}  \right]   }{\lambda_k(t_k)}\right) + \int_{t_1}^{t_k} \mu_k(s) ds  \bigg)}{f^{(1)}(z^{(1)}(1,t_1),t_1)}
\end{align}
where $f^{(1)}(z^{(1)}(1,t_1),t_1)$ is the inflow at time $t_1$ into the inner node $v_1$. Assuming no updates, we end up in a special case of (\ref{distrPa}) and set $\hat{t}_j = \hat{t}$ for the whole time period.
To improve the performance of the system it might be also possible to use more recent information at time $t_1$ for the distribution parameters (especially if $t_\text{in} < \hat{t}_j \leq t_1$) which would lead to an additional optimization problem with inflow constraint. There is an extensive study for this case under the assumption of a linear damping function in \cite{GoettlichSchillinger}. However, we remark that in the case of non-linear damping an explicit calculation is not possible anymore.

\section{Numerical Experiments}
\label{numSim}

Finally, we numerically compare the solutions to \eqref{ControlProblem11} for different non-linear damping functions $\hat{g}_i(z)$. Classical choices for such non-linear damping functions are monomials (\cite{Ikeda}), i.e. functions of the type 
\begin{align*}
	\hat{g}^{(n)}(z) = C_n z^n,
\end{align*}
where $C_n$ is a constant depending on the degree which has to be determined. By construction it holds that $\hat{g}^{(n)}(0) = 0$. For better illustration purposes, we choose the same damping functions on all three arcs. To be able to compare the damping functions we request that 
\begin{align*}
	\int_0^\frac{1}{10} \hat{g}^{(n)}(z) dz = \frac{1}{20}, ~~~ n \in \mathbb{N}
\end{align*}
and choose the constant $C_n$ accordingly. For our analysis we compare monomials up to order 4, which then are given by
\begin{align}
\begin{split}
\label{monomials}
    &\hat{g}^{(1)}(z) =  z, ~~~\hat{g}^{(2)}(z) =  15z^2, \\
	&\hat{g}^{(3)}(z) =  200z^3, ~~~\hat{g}^{(4)}(z) =  2500z^4.
\end{split}
\end{align}
Additionally we consider scenarios without any damping, i.e $\mu_i(t) = 0$.

For the demand process we choose two Jacobi processes with the parameters shown in Table \ref{paramsDN}.

\begin{table}[ht!]
	\centering
	\caption{Parameter choices for the demand processes of the 1-2 network.}
	\scalebox{0.9}{
		\begin{tabular}{c|c|c|c|c}
			\hline
			node &$\theta^{(v_i)}(t)$ & $\kappa^{(v_i)}$ & $\sigma^{(v_i)}$ & $(d_0)^{(v_i)}$    \\
			\hline
			\hline
			$v_2$& $0.45 + 0.2\sin(\pi t+1)$ & $2$& $0.9$ & $0.4$\\
			\hline
			$v_3$& $0.5 + 0.3\sin(\pi t-0.5)$ &$1$ & $0.6$ & $0.6$\\
			\hline 
	\end{tabular}}
	\label{paramsDN}
\end{table}

For the simulation of the demands we use a truncated Euler-Maruyama scheme on a time grid $(t_j)_{j \in \mathbb{N}}$ with temporal stepsize $\Delta t = \frac{1}{1000}$, i.e.
\begin{align*}
    D^*_{j+1} = D_j + \Delta t \kappa (\theta(t_j) - D_j) + \sigma \sqrt{\Delta t D_j (1-D_j)} X_j,
\end{align*}
where $X_j$ is a realization of a standard normal distributed random variable. To avoid values outside the interval $[0,1]$ for $D_{j+1}$ due to the unboundedness of $X_j$, we add a truncation into the Euler-Maruyama scheme such that the process $D_{j+1}^*$ is reflected into $[0,1]$:
\begin{align*}
    D_{j+1} &= \begin{cases} 1, &D^*_{j+1} \geq 1\\
    D_{j+1}^*, & D^*_{j+1} \in (0,1) \\
    0, & D^*_{j+1} \leq 0.\end{cases}
\end{align*}
The time horizon is chosen to $T=2.5$ and we consider an update strategy at the source node with 7 equidistant updates. 
The flux functions on the arcs are all linear in the quantity and the temporal factor of the damping function $\mu_i(t)$ obey the parameters shown in Table \ref{paramsArc}. 
\begin{table}[ht!]
	\centering
	\caption{Parameter choices for the dynamics on the arcs of the 1-2 network.}
	\scalebox{0.9}{
		\begin{tabular}{c|c|c}
			\hline
			arc & $\lambda_i(t)$ & $\mu_i(t)$ \\
			\hline
			\hline
			1 & $14 + \sin(2\pi t)$ &$1 + 0.2\sin(\pi t)$ \\
			\hline
			\hline
			2 & $12 + \sin(2\pi t)$&$0.8+ 0.2\sin(\pi t)$ \\
			\hline
			3 & $12 + \sin(4\pi t)$ &$0.8+ 0.2\sin(\pi t)$ \\
			\hline 
	\end{tabular}}
	\label{paramsArc}
\end{table}

For the discretization of \eqref{PDEwD} we use an adaptive upwind-scheme with spatial stepsize $\Delta x = \frac{1}{200}$ on a grid $(x_l)_{l \in \mathbb{N}}$. The damping is incorporated by a splitting algorithm which first performs the step from the upwind-scheme and then applies the damping in a second step. The temporal step sizes are chosen such that the CFL-condition is satisfied with equality in every time step, i.e. $\tfrac{\Delta t_j^{(i)}}{\Delta x} \lambda_i(t_j^{(i)}) = 1$. Therefore, the temporal grids depend on the velocity functions of the particular arc:
\begin{align*}
    \tilde{z}_l^{(i),j+1} = z_l^{(i),j} + \frac{\Delta t^{(i)}_j}{\Delta x} \lambda_l(t_j^{(i)}) \left(z_l^{(i),j} - z_{l-1}^{(i),j} \right).
\end{align*}
In a second step we take into account the damping and calculate $z_l^{(i),j+1}$ by
\begin{align*}
    z_l^{(i),j+1} = \tilde{z}_l^{(i),j+1} - \Delta t_j^{(i)} g_i\left(\tilde{z}_l^{(i),j+1}, t_{j+1}^{(i)}\right).
\end{align*}

\begin{figure}[ht!]
	\centering
	\includegraphics[scale = 0.55]{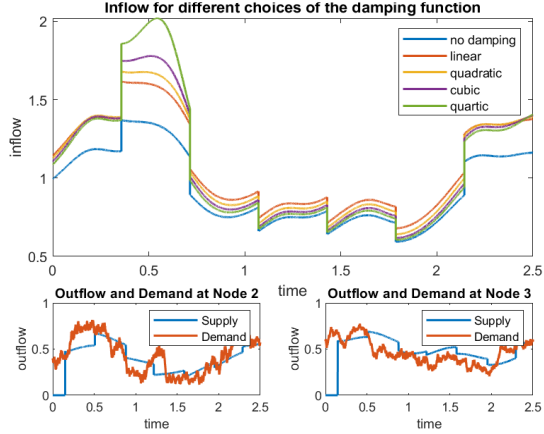}
	\caption{Comparison of different inflows for different choices of the damping function and the corresponding supplies and demands for one realization.}
	\label{dampingEx1}
\end{figure}

In Figure \ref{dampingEx1} we compare two particular representations of such a 1-2 network and show the inflows for the five different types of damping (from no damping to a quartic damping monomial in (\ref{monomials})) in the upper part and the supplies and demands at the two demand nodes for a simulation with linear damping in the lower area. Note that the supplies and demands do not differ significantly among the different damping settings. The shape of the curves is mainly influenced by the temporal non-linearities in the flux functions and the sinusoidal rhythm of the demand mean reversion levels. First, we observe that in a setting without damping the inflow is always below the inflows of the damped setting. As a second observation, we see that the inflows for the monomials are always ordered according to the degrees of $\hat{g}$. Whether they are sorted in an ascending or descending way depends on the amount of the inflow and the corresponding quantities. For larger quantities the inflows for the quartic damping function are the highest, whereas for lower inflows the inflows resulting from linear damping have the highest values.

\begin{figure}[ht!]
	\centering
	\includegraphics[scale = 0.5]{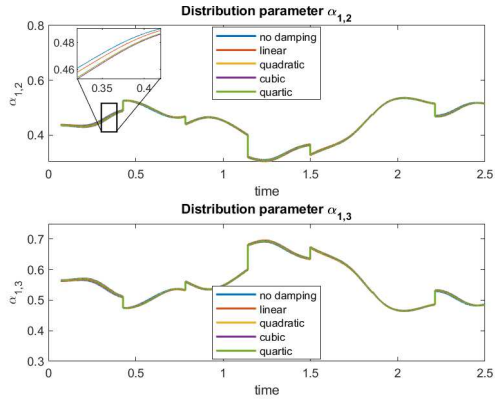}
	\caption{Comparison of the two distribution parameters for different choices of the damping function for one realization.}
	\label{dampingEx4}
\end{figure}

In Figure \ref{dampingEx4} for the same realisation that is used in Figure \ref{dampingEx1}, we provide the distribution parameters from \eqref{distrPa} at node $v_1$ for the different damping approaches. We can observe that they do not differ significantly between the different choices for the damping. In the zoomed-in frame for $\alpha_{1,2}$ we detect slight differences and that $\alpha_{1,2}$ around $t=0.35$ is a little larger for the undamped and linearly damped case.

\begin{figure}[ht!]
	\centering
	\includegraphics[scale = 0.60]{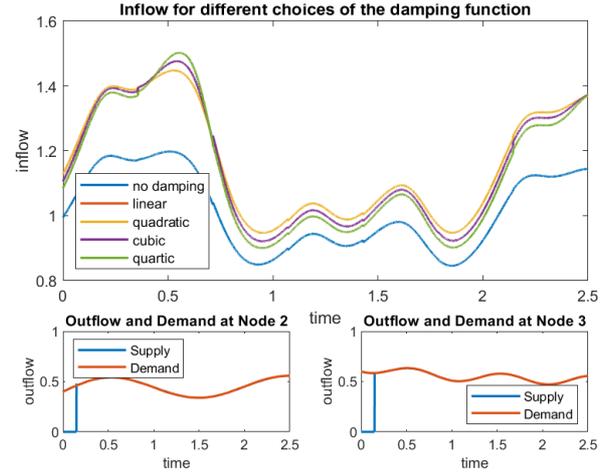}
	\caption{Comparison of different inflows for different choices of the damping function and the corresponding supplies and demands for averaged values.}
	\label{dampingEx3}
\end{figure}

In contrast to before we show the inflows in an averaged regime, where we perform a Monte Carlo simulation of $10^4$ runs in Figure \ref{dampingEx3}. First of all we observe that on average the supply matches the demand for both demand nodes. The jumps in the inflow resulting from the demand updates in the individual realisation are now almost smoothed out. But we still observe the ordering according to the damping degrees in the inflows.

\section{Conclusion}
In this work, we have presented a backward calculating procedure to determine the optimal inflow into a tree network, where the dynamics is governed by transport equations allowing for non-linear damping. In contrast to \cite{GoettlichSchillinger}, explicit calculations are more involved but still possible since the characteristic curves still do not intersect. Future work will be concerned with non-linear flux functions, where characteristic curves might cross and therefore backward calculations as done in this work are not possible any more. 


\bibliography{ifacconf}             

                                                   







\end{document}